\documentclass[11pt]{article}
\usepackage[utf8]{inputenc}
\usepackage[english]{babel}
\usepackage[svgnames]{xcolor}

\usepackage{amsthm,amsmath,amsfonts,amssymb}

\usepackage{array,enumerate,titling,float,eucal, titlesec,graphicx, textgreek, mathabx, moreverb, tikz}

\usepackage[left=2.5 cm,right=2.5cm,top=2.5cm,bottom=2.5cm]{geometry}
\usepackage[colorlinks=true,linkcolor=Purple,urlcolor=Teal,citecolor=Olive]{hyperref}

\title{Emergence of extended states at zero in the spectrum of sparse random graphs}
\author{Simon Coste \\ \url{coste@lpsm.paris} \\ \and Justin Salez \\ \url{salez@lpsm.paris}}
\date{\today}

\titleformat{\section}[block]{\filcenter}{\Large\bf\thesection .}{.5em}{\Large\bf}
\titleformat{\subsection}[block]{\filcenter}{\bf \thesubsection  .}{.5em}{\bf}
\titleformat{\paragraph}[runin]{\itshape\normalsize}{\theparagraph}{}{}
\newtheoremstyle{thm}{12pt}{7pt}{ \slshape}{}{\bfseries}{ --- }{ }{}
\theoremstyle{thm} 
\newtheorem{conjecture}{Conjecture}
\newtheorem{theorem}{Theorem}
\newtheorem{prop}[theorem]{Proposition}
\newtheorem{problem}[theorem]{Problem}
\newtheorem{question}[theorem]{Question}
\newtheorem{lem}[theorem]{Lemma}
\newtheorem{corol}[theorem]{Corollary}

\newcommand{\argmax}{\mathrm{argmax\,}}
\newcommand{\e}{\mathfrak{e}}
\newcommand{\s}{\mathfrak{s}}
\newcommand{\cGs}{{\mathcal{G}_\star}}
\newcommand{\cGss}{{\mathcal{G}_{\star\star}}}
\newcommand{\dC}{\mathbb {C}}
\newcommand{\cL}{{\mathcal{L}}}
\newcommand{\cE}{{\mathcal{E}}}

\newcommand{\cN}{{\mathcal{N}}}

\newcommand{\R}{{\mathbb{R}}}
\newcommand{\PP}{{\mathbf{P}}}
\newcommand{\EE}{{\mathbf{E}}}

\newcommand{\UGWT}{\mathsf{UGW}}

\numberwithin{equation}{section}

\begin{document}

\maketitle

\begin{abstract}
We confirm the long-standing prediction that $c=e\approx 2.718$ is the threshold for the emergence of  a non-vanishing absolutely continuous part  (extended states) at zero in the limiting spectrum of the Erd\H{o}s-Renyi random graph with average degree $c$. This is achieved by a detailed second-order analysis of  the resolvent $(A-z)^{-1}$  near the singular point $z=0$, where $A$ is the adjacency operator  of the Poisson-Galton-Watson tree with mean offspring $c$.   More generally, our method applies to arbitrary unimodular Galton-Watson trees, yielding explicit criteria for the presence or absence of extended states at zero in the limiting spectral measure of a variety of random graph models, in terms of the underlying degree distribution. 
\end{abstract}

\bibliographystyle{alpha}
\tableofcontents
\newpage
\begin{figure}
\begin{center}
\includegraphics[scale=0.5]{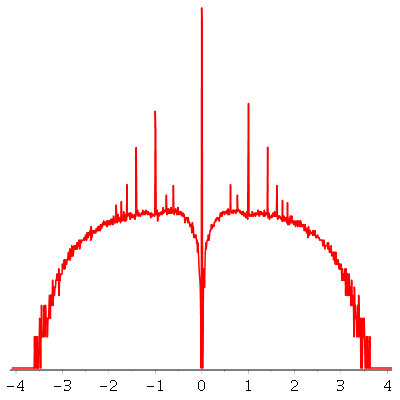}
\includegraphics[scale=0.5]{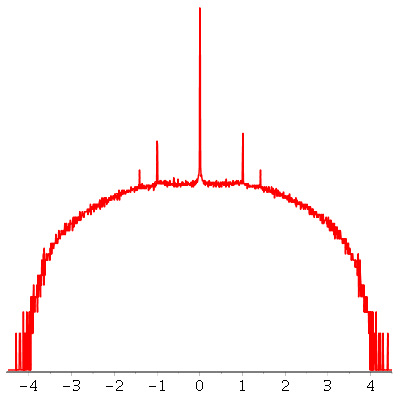}
\caption{Logarithmic plots of the adjacency spectrum of an Erd\H{o}s-Renyi random graph of size $n=10000$ and average degree $c=2$ (left) or $c=3$ (right). The presence or absence of an absolutely continuous part at zero in the $n\to\infty$ limit is already {manifest} on these finitary plots.}
\label{fig:main}
\end{center}
\end{figure}
\section{Introduction}
This paper deals with the general question of existence of a non-trivial absolutely continuous part at zero in the adjacency spectrum of unimodular Galton-Watson trees. To motivate our work, let us first briefly describe its implications for the Erd\H{o}s-Renyi random graph.
\subsection{The Bauer-Golinelli prediction}
\label{sec:ER}
Let $G_n$ be an Erd\H{o}s-Renyi random graph with size $n$ and density $p\in(0,1)$. Its adjacency matrix $A_n$ is a random symmetric $n\times n$ matrix  with zero entries along the diagonal and independent Bernoulli$(p)$ entries above the diagonal. The associated empirical eigenvalue distribution is
 \begin{eqnarray}
 \mu_{G_n} & := & \frac{1}{n}\sum_{k=1}^n\delta_{\lambda_k},
 \end{eqnarray}
where $\lambda_1\ge\ldots\ge \lambda_n$ are the eigenvalues of $A_n$. When $n\to\infty$ while $p$ is kept fixed, a celebrated result of Wigner \cite{MR1366418} asserts that a suitably rescaled version of $\mu_{G_n}$ converges weakly in probability to the semi-circle law. This remains true if $p=p_n$ tends to $0$ as $n\to\infty$, as long as $np_n\to\infty$ (see, e.g. \cite{MR2999215}). The situation changes significantly, however, when instead,
\begin{equation}
\label{sparse}
np_n\xrightarrow[n\to\infty]{} c\in(0,\infty).
\end{equation}
In  this \emph{sparse} regime, the semi-circle law gives place to a non-explicit, densely-discontinuous measure $\mu_c$, discovered in \cite{MR2259200,ksv2004} and later identified in \cite{bordenave_lelarge} as the expected {spectral measure} of the Poisson-Galton-Watson tree with mean offspring $c$ (see below). The latter has attracted a considerable attention  \cite{bordenave_lelarge_salez,bordenave2013mean,MR3315609,MR3521277,2017arXiv171007002J}, as it captures the asymptotics of many properties of $G_n$. One emblematic example is the nullity $\dim\ker(A_n)$, which is known to satisfy
\begin{equation}
\label{nullity}
\frac{1}{n}\dim\ker(A_n) \xrightarrow[n\to\infty]{\PP} \mu_c(\{0\}).
\end{equation}
In a remarkable work \cite{bauer2001exactly}, physicists Bauer and Golinelli used the so-called \emph{replica-symmetric ansatz} to predict the following intriguing formula for the limit in (\ref{nullity}).
\begin{conjecture}[Atomic mass at zero]\label{conj:rank}For any $c\in(0,\infty)$, 
\begin{equation}
\label{bg}
\mu_c(\{0\}) = q(c)+e^{-cq(c)}+cq(c)e^{-cq(c)} -1,
\end{equation}
where $q(c)$ denotes the smallest point $q\in(0,1)$ satisfying the fixed-point equation
\begin{equation}
\label{bgfp}
q=e^{-ce^{-cq}}.
\end{equation}
\end{conjecture}

\begin{figure}
\begin{center}
\includegraphics[scale=0.5]{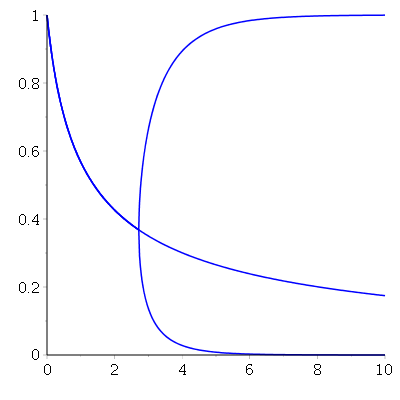}
\caption{The set of pairs $(c,q)$ satisfying the  equation (\ref{bgfp}). The branch point $(e,e^{-1})$ causes a rupture of analyticity in the spectral mass at zero $\mu_c(\{0\})$, as $c$ reaches $e$.}
\label{fig:fp}
\end{center}
\end{figure}

A quick analysis of (\ref{bgfp}) -- or an even quicker look at  Figure \ref{fig:fp} -- reveals that the right-hand side of (\ref{bg}) undergoes a  rupture of analyticity as $c$ reaches the value $e\approx 2.718$. Bauer and Golinelli proposed an interpretation of this anomaly as a phase transition in the asymptotic structure of the kernel of $A_n$. Guided by numerical simulations, they further predicted the point $c=e$ to be the threshold for the emergence of a continuous part at zero in the limiting measure $\mu_c$ \cite{bauer2001exactly,BG_incidence}. To be more precise, we will say that a measure $\mu$ has \emph{no extended states} at a location $E\in\R$ if 
\begin{eqnarray}
\frac{\mu\left([E-\varepsilon,E+\varepsilon]\right)-\mu\left(\{E\}\right)}{\varepsilon} & \xrightarrow[\varepsilon\to 0+]{} & 0,
\end{eqnarray}
and has \emph{extended states} at $E$ otherwise. This terminology is borrowed from the theory of random Schr\"odinger operators (see, e.g., \cite{MR2257129} for a recent treatment).
\begin{conjecture}[Emergence of extended states at zero]\label{conj:main} The following phase transition occurs: 
\begin{enumerate}
\item If $c<e$, then $\mu_c$ has no extended states at $0$. 
\item If $c>e$, then $\mu_c$ has extended states at $0$.
\end{enumerate}
\end{conjecture}
Conjecture \ref{conj:rank} was established almost a decade ago \cite{bordenave_lelarge_salez} by a detailed \emph{first-order} analysis of the random operator $(A-z)^{-1}$  near the singular point $z=0$, where $A$ is the adjacency operator  of the Poisson-Galton-Watson tree with mean offspring $c$.  To the best of our knowledge however, Conjecture \ref{conj:main} -- reiterated in \cite{bordenave_lelarge_salez} -- had so far remained open. In the present work, we establish this long-predicted phase transition, illustrated on Figure \ref{fig:main}. This is achieved by investigating the \emph{second-order} behavior of the random operator $(A-z)^{-1}$  near $z=0$. As already mentioned, our result is not limited to the Erd\H{o}s-Renyi model: we provide general, explicit criteria for the presence or absence of extended states at zero in the limiting spectral measure of any graph sequence whose {local weak limit} is a unimodular Galton-Watson tree, as defined next. 
\subsection{General framework}

The purpose of this section is to introduce our main objects of study, namely  spectral measures of unimodular Galton-Watson trees. We only recall the necessary notions, and refer to the comprehensive survey \cite{BorSpectrum}  for more details on graph limits and their spectral theory.
\paragraph{Spectral measures.} Let $G=(V,E)$ be a countable, locally finite graph. Its adjacency operator $A$ is a symmetric linear operator on the Hilbert space $\ell^2_{\dC}(V)$. The domain of $A$ consists of all finitely-supported vectors, and the action of $A$ on the canonical basis $(\e_x\colon x\in V)$ is given by
\begin{equation}
\langle  \e_x | A\, \e_y\rangle =
\left\{
\begin{array}{ll}
1 & \textrm{if }\{x,y\}\in E\\
0 & \textrm{otherwise.}
\end{array}
\right.
\end{equation}
As long as $A$ is essentially self-adjoint, the Spectral Theorem applies: the resolvent $(A-z)^{-1}$ is a well-defined bounded operator for all $z\in\dC\setminus\R$, and for every $o\in V$, we have the representation
\begin{equation}
\label{def:rootedmeasure}
\forall z\in \dC\setminus\R,\qquad \langle \e_o|(A-z)^{-1}\e_o \rangle = \int_\R\frac{1}{\lambda-z}\,\mu_{(G,o)}(\mathrm{d}\lambda),
\end{equation}
for a unique probability measure  $\mu_{(G,o)}$ on $\R$, called the spectral measure of the rooted graph $(G,o)$. This fundamental object will be central to our work. It may be thought of as the local contribution of  $o$ to the spectrum of $G$. Indeed, when $G$ is finite, there is an orthonormal basis of $n=|V|$ eigenfunctions $\phi_1,\ldots,\phi_{n}$ of $A$ with respective eigenvalues $\lambda_1,\ldots,\lambda_{n}$, and we have the   expression
\begin{equation}
\mu_{(G,o)}=\sum_{k=1}^{n}|\phi_k(o)|^2\delta_{\lambda_k}.
\end{equation}
In particular, the empirical eigenvalue distribution  $\mu_G:=\frac 1n\sum_{k=1}^n\delta_{\lambda_k}$ can be recovered from the spectral measures $(\mu_{(G,o)}\colon o\in V)$ by averaging over the choice of the root: 
\begin{equation}
\mu_G = \frac{1}{|V|}\sum_{o\in V}\mu_{(G,o)}.
\end{equation}
Of course, neither side of this identity makes sense when $G$ is infinite. However, the framework of \emph{local weak convergence} enables us  to pass to the ``infinite-volume limit'', in an appropriate sense.

\paragraph{Local weak convergence.}   Write $\cGs$ for the space of locally finite, connected rooted graphs, considered up to root-preserving isomorphism. Make it complete and separable  by letting the distance between two rooted graphs be $1/(1+r)$, where $r$ is the largest integer such that the balls of radius $r$ around the root in the two graphs are isomorphic. Any finite graph naturally induces a probability measure on $\cGs$, via choosing a root uniformly at random and restricting to its connected component. If $(G_n)_{n\ge 1}$ is a sequence of finite graphs, and if the sequence of probability measures thus induced admits a weak limit $\cL$, then $\cL$ is called the \emph{local weak limit} of  $(G_n)_{n\ge 1}$ \cite{aldous_steele,benjamini_schramm}. In words, $\cL$ is the law of a random rooted graph $(G,o)$ that describes what $G_n$ asymptotically looks like when seen from a uniformly chosen vertex.

This limiting object has been shown to capture a number of asymptotic properties of $(G_n)_{n\ge 1}$, including the empirical eigenvalue distribution itself. More precisely, we always have
\begin{equation}
\label{lwc}
\sup_{\lambda\in\R}\left|\mu_{G_n}\left((-\infty,\lambda]\right)-\mu_\cL\left((-\infty,\lambda]\right)\right|\xrightarrow[n\to\infty]{} 0,
\end{equation}
where $\mu_{\cL}(\cdot):=\EE\left[\mu_{(G,o)}(\cdot)\right]$ denotes the expected spectral measure under $\cL$\footnote{This definition implicitly relies on the (non-trivial) fact that the adjacency operator of a unimodular random graph is essentially self-adjoint with probability $1$,  see \cite[Proposition 2.2]{BorSpectrum} for a proof.}.
 This remarkable continuity principle has a long history  \cite{bordenave_lelarge,pointwise,BorSpectrum}. In short, it allows one to replace the spectral analysis of sparse graphs by that of their local weak limits. Luckily, the latter turn out to be much more convenient to work with than the finite graphs that they approximate. For example, although they  have many cycles, most sparse random graphs admit a local weak limit that is supported on trees. Moreover, in many cases of interest, including the Erd\H{o}s-Renyi and configuration models, the limit has a particularly simple recursive structure, which we now describe.

\paragraph{Unimodular Galton-Watson trees.}
Let $\pi=(\pi_k)_{k\geq 0}$ be a probability distribution on $\mathbb N$ with finite, non-zero  mean. 
A {unimodular Galton-Watson tree} with degree distribution $\pi$ is a random rooted tree obtained by a Galton-Watson branching process in which the root has offspring distribution $\pi$ and all descendants have the size-biased offspring distribution  $\widehat{\pi}=(\widehat{\pi}_k)_{k\geq 0}$ given by
\begin{equation}
\label{eq:sizebiased}
\widehat\pi_k := \frac{(k+1)\pi_{k+1}}{\sum_{i}i\pi_{i}}.
\end{equation}
The law of this random rooted tree plays a distinguished role in the theory and will be denoted by $\UGWT(\pi)$. It arises as the local weak limit of uniform random graphs with prescribed degrees, when the number of vertices tends to infinity while the empirical degree distribution tends to $\pi$. 

\begin{figure}
\begin{center}
\includegraphics[angle =0,height =6.8cm]{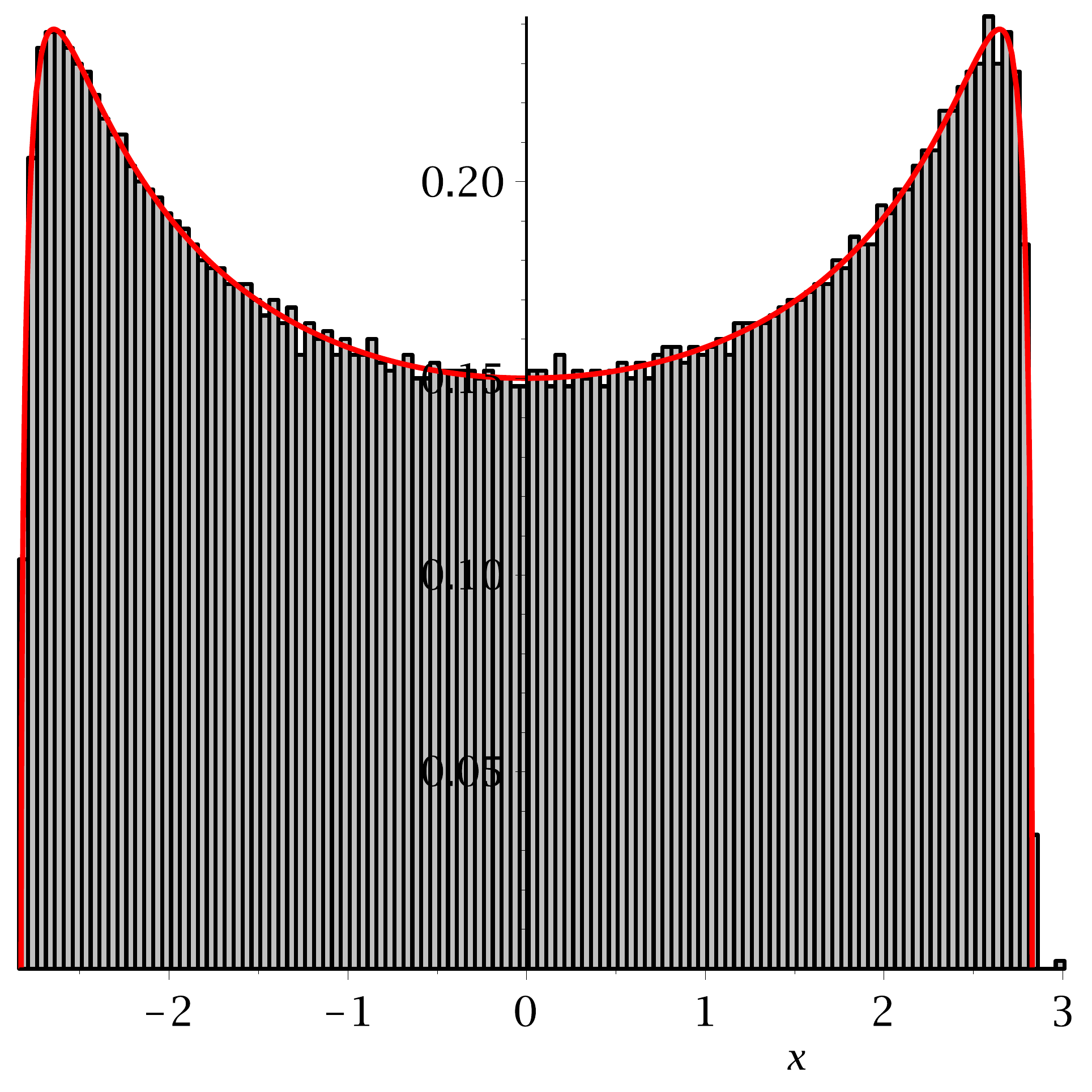}
\caption{Histogram of the eigenvalues of a uniform $3-$regular random graph on $10^4$ vertices (in gray), and the limiting Kesten McKay density (in red)\label{fig:mckay}.}
\end{center}
\end{figure}

A simple example is random $d-$regular graphs, for which $\pi$ is just a Dirac mass at $d$: the resulting tree is then the infinite $d-$regular rooted tree, whose spectral measure is the well-known Kesten-McKay distribution \cite{mckay}, see Figure \ref{fig:mckay}. Another important example is the Erd\H{o}s-Renyi model with parameters as in (\ref{sparse}), for which $\pi$ is the Poisson distribution with mean $c$. In that case, we have $\widehat{\pi}=\pi$, so that $\UGWT(\pi)$ is the law of the standard Poisson-Galton-Watson tree with  mean offspring $c$. Its expected spectral measure $\mu_{\UGWT(\pi)}$ is precisely the limit $\mu_c$ mentioned in Section \ref{sec:ER}. The striking difference in the spectra of these two models (see Figures \ref{fig:main} and \ref{fig:mckay}) motivates the following research program, to which the present paper is intended to contribute.
\begin{problem}\label{problem}Understand the regularity of $\mu_{\UGWT(\pi)}$ -- in particular, the  supports of its pure-point, absolutely-continuous, and singular-continuous parts -- as a function of the degree distribution $\pi$. 
\end{problem}

\paragraph{State of the art.} 
This relatively young line of research has already witnessed notable progress. A comprehensive account, as well as a list of exciting conjectures, can be found in the introductory survey \cite{BorSpectrum}.  The pure-point part of the spectrum is now reasonably well understood. In particular, the work \cite{bordenave_lelarge_salez} provides an explicit formula  for the mass at zero, while \cite{MR3315609, spectral_atoms} investigate the locations of other atoms. Rigorous results on the support of the continuous part are more limited. A remarkably general criterion by Bordenave, Sen and Vir\'ag \cite{bordenave2013mean} guarantees the presence of a non-trivial continuous part as soon as the Galton-Watson tree is super-critical. Unfortunately, the result is existential in nature and can not be used to ensure the presence of extended states at a given location. More precise information is available when $\pi$ is sufficiently close to a Dirac mass, thanks to Keller \cite{MR2994759}. However, the method used there is intrinsically perturbative and does not yield information for explicit choices of $\pi$ such as the Poisson distribution involved in Conjecture \ref{conj:main}. In the present paper, we provide explicit criteria for the presence or absence of extended states at zero in the spectral measure $\mu_{\UGWT(\pi)}$, for a general degree distribution $\pi$.

\subsection{Results}
\label{sec:main}
Throughout this section, we fix a probability measure $\pi$ on $\mathbb N$ with finite, non-zero mean, and we let  $\mu=\mu_{\UGWT(\pi)}$ denote the expected spectral measure of the unimodular Galton-Watson tree with degree distribution $\pi$. In the degenerate case where $\pi_0+\pi_1=1$, our random tree is just an isolated vertex with probability $\pi_0$ and an isolated edge with probability $\pi_1$, so its expected spectral measure is $\mu_{\UGWT(\pi)}=\pi_0\delta_0+\frac{\pi_1}{2}\delta_{-1}+\frac{\pi_1}{2}\delta_{1}$, which trivially has no extended states anywhere. To avoid degeneracies, we will henceforth always assume that 
\begin{align}
\label{degen}
\pi_0+\pi_1<1.
\end{align}
All our results will be expressed in terms of the degree generating series
\begin{equation}
\varphi(z) := \sum_{k\ge 0}\pi_k z^k,\qquad \widehat{\varphi}(z):=\sum_{k\ge 0}\widehat{\pi}_k z^k = \frac{\varphi'(z)}{\varphi'(1)}.
\end{equation}
It was shown in \cite[Theorem 2]{bordenave_lelarge_salez} that $\mu(\{0\}) = \max M$, where the function $M\colon [0,1]\to \R$ is given by
\begin{equation}
M(z) := \varphi(z)+(1-z)\varphi'(z)+\varphi\left(1-\widehat{\varphi}(z)\right)-1.
\end{equation}
Our main finding is that the presence or absence of extended states at zero depends on the set 
\begin{equation}
\argmax M := \left\{z\in[0,1]\colon M(z)=\max M\right\}.
\end{equation}
A quick differentiation shows that any $z\in\argmax M$ must satisfy 
\begin{equation}
\label{fp}
z = 1-\widehat{\varphi}\left(1-\widehat{\varphi}(z)\right).
\end{equation}
Among the (possibly many) solutions to this fixed-point equation, the following one will play a crucial role: we let $z_\star\in(0,1)$ denote the unique point satisfying
\begin{equation}
\label{zstar}
z_\star = 1-\widehat{\varphi}\left(z_\star\right).
\end{equation}
It is easily checked that $M'(z_\star)=0$, and that $M''(z_\star)$ has the same sign as $\widehat{\varphi}'(z_\star)-1$. 
The presence or absence of extended states at zero in $\mu$ turns out to be dictated by the following two conditions:
 \begin{enumerate}[(i)]
 \item $M$ achieves its maximum uniquely at $z_\star$, i.e. $\argmax M=\{z_\star\}$.
 \item $M''(z_\star)\ne 0$ (or equivalently, $\widehat{\varphi}'(z_\star)\ne 1$).
 \end{enumerate}
More precisely, our first main result states that (i) and (ii) characterize a strong square-integrability property which, in particular, implies the absence of extended states at zero. 
\begin{theorem}[No extended states at zero]\label{th:noext}The square-integrability property
\begin{equation}
\label{squareint}
\int_{\R\setminus\{0\}}\frac{1}{\lambda^2}\,\mu(\mathrm{d}\lambda)<\infty
\end{equation}
 holds if and only if (i) and (ii) are both satisfied. 
In particular, when this is the case, $\mu$ satisfies
\begin{equation}
\label{noextstrong}
\mu\left(\left[-\varepsilon,\varepsilon\right]\right)=\mu\left(\{0\}\right) + o(\varepsilon^2),
\end{equation}
as $\varepsilon\to 0$, which is much stronger than the absence of extended states at zero.
\end{theorem}
Conversely, our second main result guarantees the existence of extended states at zero as soon as (i) fails. Note that this only leaves aside the critical situation where (i) holds but (ii) fails, in which case we do not know whether the measure $\mu$ has extended states at zero. We emphasize that this situation is not generic, as it forces  $z\mapsto \widehat{\varphi}(z)+z$ and $z\mapsto \widehat{\varphi}'(z)$ to reach $1$ at the same point.
\begin{theorem}[Extended states at zero]\label{th:ext}If condition (i) fails, then $\mu$ has extended states at zero.
\end{theorem}

We end this section by applying our results to the special case where $\pi$ is the Poisson distribution with mean $c$, i.e. $\varphi(z)=\widehat{\varphi}(z)=e^{c(z-1)}$. Under the change of variable $q=1-z$, the fixed-point equation (\ref{fp}) reduces to (\ref{bgfp}), whose solutions were represented on Figure \ref{fig:fp}. When $c\le e$, the solution is unique, so condition (i) trivially holds. When $c> e$, there are three solutions, and condition (i) must fail because $M''(z_\star)>0$. In fact, the double equality $\widehat{\varphi}'(z_\star)=ce^{c(z_\star-1)}=c(1-z_\star)$ shows that $M''(z_\star)$ is negative, null or positive according to whether $c$ is less than, equal to, or more than $e$. Thus, Theorem \ref{th:noext} applies if and only if $c<e$, and Theorem \ref{th:ext} applies if and only if $c>e$. This establishes Conjecture \ref{conj:main} and leaves aside the critical case $c=e$, which remains open.

\iffalse
\begin{question}In the critical case where the mean degree is $c=e$, does the expected spectral measure of the Poisson-Galton-Watson tree has extended states at zero?
\end{question}

%

Bauer and Golinelli : \cite{BG_kernel_tree}

Spectrum of random trees : \cite{spectrum_tree}

\fi

\section{Main ingredients}
\label{sec:prelim}
In this section, we introduce the main ingredients of our proof. We start by reformulating the problem of extended states at zero in terms of Stieltjes transforms, and then recall the well-known recursion satisfied by the latter on rooted trees. We then combine this recursion with the Mass Transport Principle to establish a new identity that will be crucial to our proof. We emphasize that all results in this section apply to general trees. The special structure of unimodular Galton-Watson trees will only enter the play in Section \ref{sec:last} below.
\subsection{Stieltjes transform}\label{section:def}
The integral appearing in the definition (\ref{def:rootedmeasure}) is known as the \emph{Stieltjes transform} of the measure $\mu_{(G,o)}$. Here we will focus on the imaginary part of its restriction to the imaginary axis. More precisely, given a finite Borel measure $\mu$ on $\R$, we consider the observable $\s\colon (0,\infty)\to\R$ defined by
\begin{equation}
\label{stieltjes}
\s(t) :=  \int_\R\frac{t}{\lambda^2+t^2}\,\mu(\mathrm{d}\lambda).
\end{equation}
The leading-order asymptotics of $\s(t)$ as $t\to 0$ are directly related to the behavior of $\mu$ around $0$. In particular, the following two limits emerge naturally:
\begin{equation}
\alpha:=\lim_{t\to 0}\downarrow t\s(t)=\mu\left(\{0\}\right), \qquad
\beta:=  \lim_{t\to 0}\uparrow\frac{\s(t)}{t} =\int_\R\frac{1}{\lambda^2}\,\mu(\mathrm{d}\lambda).
\end{equation}
Note that we can not simultaneously have $\alpha>0$ and $\beta<\infty$. To investigate the presence or absence of extended states at zero, we first need to subtract the atom at zero from $\mu$, i.e. consider the measure $\mu^\star:=\mu-\alpha \delta_0$ and its associated transform, 
\begin{equation}
\s^\star(t)=\s(t)-\frac{\alpha}{t}=\int_{\R\setminus\{0\}}\frac{t}{\lambda^2+t^2}\,\mu(\mathrm{d}\lambda).
\end{equation}
We then have the following exact characterization of the absence of extended states at zero.
\begin{lem}[Characterization]\label{lm:ext}$\mu$ has no extended states at zero if and only if $\s^\star(t)\to 0$ as $t\to 0$.
\end{lem}
\begin{proof} Fix $t>0$ and set $I_t:=[-t,t]\setminus\{0\}$. Since $\frac{t}{\lambda^2+t^2}\ge \frac{1}{2t}$ for all $\lambda\in I_t$, we have
\begin{align*}
\s^\star(t) &\ge  \frac{\mu\left(I_t\right)}{2t}.
\end{align*}
Thus, $\mu$ has no extended states at zero whenever $s^\star(t)\to 0$. Conversely, observe that for any $\varepsilon,t>0$, we have by Fubini's Theorem
\begin{align*}
\int_{I_\varepsilon}\frac{t}{\lambda^2+t^2}\,\mu(\mathrm{d}\lambda)
&=\int_0^\infty \frac{2tu}{(t^2+u^2)^2}\,\mu\left(I_{\varepsilon\wedge u}\right)\mathrm{d}u.
\end{align*}
On the other hand, the same identity with the measure $\mu$ replaced by Lebesgue's measure gives
\begin{align*}
\arctan\left(\frac{\varepsilon}{t}\right) &=\int_0^\infty \frac{2tu}{(t^2+u^2)^2}\left(\varepsilon\wedge u\right)\mathrm{d}u.
\end{align*}
Comparing these two lines, we deduce that
\begin{align*}
\int_{I_\varepsilon}\frac{t}{\lambda^2+t^2}\,\mu(\mathrm{d}\lambda) &\le \sup_{u\in(0,\varepsilon)}\left\{\frac{\mu(I_u)}{u}\right\}\arctan\left(\frac{\varepsilon}{t}\right).
\end{align*}
Since $\arctan\left(\cdot\right)\le \frac{\pi}{2}$, and since $\frac{t}{\lambda^2+t^2}\le \frac{t}{\varepsilon^2}$ for all $\lambda\notin I_\varepsilon$, we conclude that
\begin{align*}
\s^\star(t) &\le \frac{t}{\varepsilon^2}+\frac{\pi}{2}\sup_{u\in(0,\varepsilon)}\left\{\frac{\mu(I_u)}{u}\right\}.
\end{align*}
Sending $t\to 0$ and then $\varepsilon\to 0$ shows that $s^\star(t)\to 0$ whenever $\mu$ has no extended states at zero.
\end{proof}
This lemma reduces the absence of extended states at zero to the condition $\s^\star(t)=o(1)$ as $t\to 0$.
Moreover, the square-integrability property (\ref{squareint}) can be rephrased as $\beta^\star<\infty$, where
\begin{equation}
\label{betastar}
\beta^\star:=  \lim_{t\to 0}\uparrow\frac{\s^\star(t)}{t} =\int_{\R\setminus\{0\}}\frac{1}{\lambda^2}\,\mu(\mathrm{d}\lambda).
\end{equation}
Thus, our two main theorems will follow from a careful analysis of $\s^\star(t)$ as $t\to 0$, when $\mu$ is the expected spectral measure of a unimodular Galton-Watson tree. The starting point of this analysis is a well-known local recursion satisfied by spectral measures of rooted trees.

\subsection{Local recursion}
As many graph-theoretical quantities, spectral measures 
admit a recursive structure when evaluated on trees.  
Fix a tree $T=(V,E)$ whose adjacency operator is self-adjoint, and let $o\in V$ be an arbitrary vertex. We write $\partial o=\{x\in V\colon \{x,o\}\in E\}$ for the set of its neighbours, and $\deg(o)=|\partial o|$ for its degree. Deleting $o$ splits $T$ into $\deg(o)$ disjoint subtrees which will naturally be denoted by $\left(T_{x\to o}\colon x\in\partial o\right)$. We let $\s_{o},\s^\star_{o},\alpha_{o},\beta_{o},\beta^\star_{o}$ be the objects $\s,\s^\star,\alpha,\beta,\beta^\star$ defined above when the general measure $\mu$ is taken to be the spectral measure $\mu_{(T,o)}$. Similarly, we let $\s_{x\to o},\s^\star_{x\to o},\alpha_{x\to o},\beta_{x\to o},\beta^\star_{x\to o}$ correspond to the choice $\mu=\mu_{(T_{x\to o},x)}$.  We then have the following elementary but fundamental relation (see, e.g., \cite{bordenave_lelarge_salez} for a proof): for all $t\in(0,\infty)$,
\begin{equation}
\label{decide}
\s_o(t)  =  \frac{1}{t+\sum_{x\in\partial o}\s_{x\to o}(t)}.
\end{equation} 
In particular, multiplying or dividing both sides by $t$ and sending $t\to 0$ yields
\begin{align}
\label{decide:alpha}
\alpha_o &=  \frac{1}{1+\sum_{x\in\partial o}\beta_{x\to o}};\\
\label{decide:beta}
\beta_o &=  \frac{1}{\sum_{x\in\partial o}\alpha_{x\to o}}.
\end{align}
In view of these identities, it is natural to decompose the degree as $\deg(o)=\cN_o^++\cN_o^-+\cN_o^\star$ where
\begin{align}
\cN_o^+:=\sum_{x\in\partial o}{\bf 1}_{(\alpha_{x\to o}>0)},\qquad \cN_o^-:=\sum_{x\in\partial o}{\bf 1}_{(\beta_{x\to o}<\infty)},\qquad\cN_o^\star:=\sum_{x\in\partial o}{\bf 1}_{(\alpha_{x\to o}=0,\beta_{x\to o}=\infty)}.
\end{align}
It then readily follows from (\ref{decide:alpha}) and (\ref{decide:beta}) that 
\begin{align}
\label{char:plus}
\alpha_o>0 &\Longleftrightarrow \cN_o^+=\cN_o^\star=0\\
\label{char:minus}
\beta_o<\infty &\Longleftrightarrow \cN_o^+\ge 1\\
\label{char:star}
\left(\alpha_o=0,\beta_o=\infty\right) &\Longleftrightarrow \left(\cN_o^+=0,\cN_o^\star\ge 1\right).
\end{align}
Of course, the recursion (\ref{decide}) also applies to the tree $T_{o\to y}$ (for any $y\in\partial o$), yielding 
\begin{align}
\label{update}
\s_{o\to y}(t)  &=  \frac{1}{t+\sum_{x\in\partial o\setminus\{y\}}\s_{x\to o}(t)};\\
\label{update:alpha}
\alpha_{o\to y} &=  \frac{1}{1+\sum_{x\in\partial o\setminus\{y\}}\beta_{x\to o}};\\ 
\label{update:beta}
\beta_{o\to y} &=  \frac{1}{\sum_{x\in\partial o\setminus\{y\}}\alpha_{x\to o}}.
\end{align}
 These recursions will play a crucial role in our analysis.
\subsection{Mass Transport Principle}
The second-order quantity $\beta^\star_o$ is \emph{a priori} much harder to analyze than its first-order counterpart $\beta_o$, as we have to remove the singularity caused by the atom at zero.  To overcome this difficulty, we will exploit a powerful identity known as the \emph{Mass Transport Principle} (see, e.g., \cite{MR2354165}): any random rooted graph $(G,o)$ whose law is the local weak limit of some sequence of finite graphs is \emph{unimodular}, in the sense that it satisfies the distributional symmetry  
\begin{equation}
\EE\left[\sum_{x\in V(G)}f(G,o,x)\right] = \EE\left[\sum_{x\in V(G)}f(G,x,o)\right],
\end{equation}
for any Borel-measurable function $f\colon\cGss\to[0,\infty]$, where $\cGss$ denotes the natural analogue of $\cGs$ for \emph{doubly-rooted} graphs. At an intuitive level, this identity expresses the fact that the root is ``equally likely'' to be any vertex (even though the underlying graph is possibly infinite). Here we use this spatial stationarity to prove the following key formula, which expresses $\EE[\beta_o^\star]$ in terms of $\beta_o$ only. 
\begin{prop}[Getting rid of the atom at zero]\label{pr:uni} For any unimodular random tree $(T,o)$, 
\begin{equation}\label{for:beta*}
\mathbf{E}[\beta^\star_o] = \mathbf{E}\left[\mathbf{1}_{(\alpha_o=0)}\beta_o \right]+\mathbf{E}\left[\mathbf{1}_{(\mathcal{N}_o^+ \geqslant 2)} \beta_o\right]+
\EE\left[\mathbf{1}_{(\mathcal{N}_o^+ \geqslant 2)}\frac{\sum_{x \in \partial o} \beta_{x \to o}\mathbf{1}_{(\alpha_{x \to o}=0)}}{\sum_{x\in\partial o}\alpha_{x \to o}\mathbf{1}_{(\alpha_{x \to o}>0)}}\right].
\end{equation}
\end{prop}
%We recall that $\{\cN_o^+\ge 2\}\subseteq \{\beta_o<\infty\}\subseteq\{\alpha_o=0\}$. Thus, there is no atom at zero inside the expectation on the right-hand side, making the identity extremely useful. 
\begin{proof}Fix $t\in(0,\infty)$ and $y\in\partial o$. Combining (\ref{decide}) and (\ref{update}), we have
\begin{equation}
{\s_o(t)}=\left(\frac{1}{\s_{o\to y}(t)} + \s_{y\to o}(t)\right)^{-1}.
\end{equation}
Multiplying by $\s_{y\to o}(t)$ clearly makes the right-hand side symmetric in $o$ and $y$, and hence
\begin{equation}
{\s_o(t)\s_{y\to o}(t)}={\s_y(t)\s_{o\to y}(t)}.
\end{equation}
Summing over all $y\in\partial o$ and using again (\ref{decide}), we obtain
\begin{align}
\label{first}
1-t\s_o(t)=\sum_{y\in\partial o}\s_y(t)\s_{o\to y}(t).
\end{align}
On the other hand, it easily follows from (\ref{decide:alpha}),(\ref{update:alpha}) and (\ref{update:beta}) that for any $y\in\partial o$, we have
\begin{equation}
\label{second}
\alpha_{o}>0 \Longleftrightarrow \left(\alpha_{o\to y}>0\textrm{ and } \cN_y^+\ge 2\right).
\end{equation} 
Combining this with (\ref{first}), we deduce that
\begin{align}
\left(1-t\s_o(t)\right){\bf 1}_{(\alpha_o>0)}=\sum_{y\in\partial o}\s_y(t)\s_{o\to y}(t){\bf 1}_{(\alpha_{o\to y}>0)}{\bf 1}_{(\cN_y^+\ge 2)}.
\end{align}
We may now take expectation and use unimodularity to obtain
\begin{align}
\label{key}
\EE\left[\left(1-t\s_o(t)\right){\bf 1}_{(\alpha_o>0)}\right]&=\EE\left[{\bf 1}_{(\cN_o^+\ge 2)}\s_o(t)\sum_{y\in\partial o}\s_{y\to o}(t){\bf 1}_{(\alpha_{y\to o}>0)}\right].
\end{align}
Letting $t\to 0$ and using $\{\cN_o^+\ge 2\}\subseteq\{\beta_o<\infty\}$ for the right-hand side, we obtain
\begin{align*}
\EE\left[\left(1-\alpha_o\right){\bf 1}_{(\alpha_o>0)}\right]&=\EE\left[{\bf 1}_{(\cN_o^+\ge 2)}\beta_o\sum_{y\in\partial o}\alpha_{y\to o}\right]\\
&=\PP\left(\cN_o^+\ge 2\right)\\
&=\EE\left[{\bf 1}_{(\cN_o^+\ge 2)}\s_o(t)\left(t+\sum_{y\in\partial o}\s_{y\to o}(t)\right)\right],
\end{align*}
where the second line follows from (\ref{decide:beta}) and the third from (\ref{decide}).
Substracting (\ref{key}), we arrive at
 \begin{align}
\EE\left[\left(t\s_o(t)-\alpha_o\right){\bf 1}_{(\alpha_o>0)}\right]&=\EE\left[{\bf 1}_{(\cN_o^+\ge 2)}\s_o(t)\left(t+\sum_{y\in\partial o}\s_{y\to o}(t){\bf 1}_{(\alpha_{y\to o}=0)}\right)\right].
\end{align}
Dividing through by $t^2$ and sending $t\to 0$ yields 
\begin{equation}
\mathbf{E}[\beta^\star_o{\bf 1}_{(\alpha_o>0)}] = \mathbf{E}\left[\mathbf{1}_{(\mathcal{N}_o^+ \geqslant 2)} \beta_o \left(1+ \sum_{x \in \partial o} \beta_{x \to o}\mathbf{1}_{\left(\alpha_{x \to o}=0\right)}\right) \right].
\end{equation}
On the other hand, on the event $\{\alpha_o=0\}$, we have $\beta_o^\star=\beta_o$, which concludes the proof.
\end{proof}

\section{Unimodular Galton-Watson trees}\label{sec:NES}
\label{sec:last}
The above results were valid for any unimodular random tree. We now consider the special case of unimodular Galton-Watson trees, and exploit their self-similar nature to turn the above recursions into distributional fixed-point equations that will be amenable to analysis. 
\subsection{Distributional fixed-point equations}
From now on, we fix a degree distribution $\pi$ as in Section \ref{sec:main}, and we equip the space of rooted trees $(T,o)$ with two different probability measures: we reserve the letter $\PP$ for the unimodular law $\UGWT(\pi)$, and use $\widehat{\PP}$ to denote the homogeneous Galton-Watson law with offspring distribution $\widehat{\pi}$. We naturally use $\EE$ and $\widehat{\EE}$ to denote the corresponding expectations. Thus, the distribution of the root-degree $\deg(o)$ is $\pi$ under $\PP$ and $\widehat{\pi}$ under $\widehat{\PP}$ and in both cases, conditionally on $\deg(o)$, the subtrees $(T_{x\to o},x\in\partial o)$ are  i.i.d. homogeneous Galton-Watson trees with offspring distribution $\widehat{\pi}$. In particular, the recursion (\ref{decide}) takes the following simple distributional form.
\begin{corol}[Distributional structure]\label{corol}Under both $\PP$ and $\widehat{\PP}$, the conditional law of $(\cN_o^+,\cN_o^-,\cN_o^\star)$ given $\deg(o)$  is Multinomial  with parameters $\deg(o)$ and $\left(\widehat{\PP}(\alpha_o>0),\widehat{\PP}(\beta_o<\infty),\widehat{\PP}(\alpha_o=0,\beta_o=\infty)\right)$. Moreover, conditionally on $(\cN_o^+,\cN_o^-,\cN_o^\star)$, the random sums 
\begin{align*}
\sum_{x\in\partial o}\s_{x\to o}(t){\bf 1}_{(\alpha_{x\to o}>0)}, \qquad \sum_{x\in\partial o}\s_{x\to o}(t){\bf 1}_{(\alpha_{x\to o}<\infty)},\qquad \sum_{x\in\partial o}\s_{x\to o}(t){\bf 1}_{(\alpha_{x\to o}=0,\beta_{x\to o}=\infty)}
\end{align*}
are independent, the first (resp. second, resp. third) being distributed as a sum of $\cN_o^+$ (resp. $\cN_o^-$, resp. $\cN_o^\star$) i.i.d. random variables with law $\widehat{\PP}\left(\s_o(t)\in\cdot|\alpha_o>0\right)$ (resp. $\widehat{\PP}\left(\s_o(t)\in\cdot|\beta_o<\infty\right)$, resp. $\widehat{\PP}\left(\s_o(t)\in\cdot|\alpha_o=0,\beta_o=\infty\right)$).
\end{corol}
We shall use this fact (and its $t\to 0$ counterparts) repeatedly below, without notice. For example, an immediate consequence of this and (\ref{char:plus})-(\ref{char:minus}) is that
\begin{align}
\label{rec:pq}
{\PP}(\alpha_o>0)={\varphi}\left( \widehat{\PP}(\beta_o<\infty)\right),\qquad 
{\PP}(\beta_o=\infty)={\varphi}\left( \widehat{\PP}(\alpha_o=0)\right),\\
\label{rec:pqhat}
\widehat{\PP}(\alpha_o>0)=\widehat{\varphi}\left( \widehat{\PP}(\beta_o<\infty)\right),\qquad 
\widehat{\PP}(\beta_o=\infty)=\widehat{\varphi}\left( \widehat{\PP}(\alpha_o=0)\right),
\end{align}
where we recall that $\varphi$ and $\widehat{\varphi}$ are the generating series of $\pi$ and $\widehat{\pi}$ respectively. 
In particular, the number $z=\widehat{\PP}(\alpha_o=0)$ must solve the fixed-point equation (\ref{fp}). In fact, $\widehat{\PP}(\alpha_o=0)$ was shown in  \cite{bordenave_lelarge_salez} to coincide with the last point at which the function $M$ achieves its maximum, i.e.
\begin{align}
\label{charac}
\widehat{\PP}(\alpha_o=0)=\max\left(\argmax M\right).
\end{align}
With this characterization in hands, we may reformulate our main assumption (i) as follows.
\begin{lem}[Reformulation of assumption (i)]\label{lm:i}The following conditions are equivalent.
\begin{enumerate}
\item $\argmax M=\{z_\star\}$;
\item $\widehat{\mathbf{P}}(\alpha_o=0)=z_\star$;
\item $\PP\left(\alpha_o=0,\beta_o=\infty\right)=0$;
\item $\widehat{\PP}\left(\alpha_o=0,\beta_o=\infty\right)=0$.
\end{enumerate}
\end{lem}
\begin{proof}
Set $z:=\widehat{\PP}(\alpha_o=0)$. Since $\{\alpha_o>0\}\subseteq\{\beta_o=\infty\}$, we always have 
\begin{align*}
\PP\left(\alpha_o=0,\beta_o=\infty\right)&=\PP\left(\beta_o=\infty\right)-\PP\left(\alpha_o>0\right)\\
&= \varphi\left(z\right)-\varphi\left(1-\widehat{\varphi}(z)\right).
\end{align*}
where the second line follows from (\ref{rec:pq}). Similarly, 
\begin{align*}
\widehat{\PP}\left(\alpha_o=0,\beta_o=\infty\right)&=\widehat{\varphi}\left(z\right)-\widehat{\varphi}\left(1-\widehat{\varphi}(z)\right).
\end{align*}
From these equalities and the fact that $\varphi,\widehat{\varphi}$ are increasing, we immediately deduce that the conditions (2),(3) and (4) are equivalent. Moreover, it is clear from (\ref{charac}) that (1) implies (2). To see that (2) implies (1),  recall that any point $z\in\argmax M$ must satisfy the fixed-point equation (\ref{fp}), and that the latter implies  $M(1-\widehat{\varphi}(z))=M(z)$. Thus, the set $\argmax M$ is stable under the map $z\mapsto 1-\widehat{\varphi}(z)$, and so it can not intersect $(z_\star,1]$ without also intersecting $[0,z_\star)$.
\end{proof}

\subsection{Proof of Theorem \ref{th:noext}}

In this section, we prove Theorem \ref{th:noext}, namely, that (i) and (ii) are  necessary and sufficient for
$
\EE[\beta_o^\star]<\infty.
$
The necessity of (i) is easy: if (i) fails, then $\PP\left(\alpha_o=0,\beta_o=\infty\right)>0$ by Lemma \ref{lm:i}, and so the first term on the right-hand side of (\ref{for:beta*}) is already infinite. We will thus henceforth assume that (i) holds. By Lemma \ref{lm:i}, this ensures that
\begin{align}
\label{now}
\widehat{\PP}(\alpha_o=0)=\widehat{\PP}(\beta_o<\infty)&=z_\star.
\end{align}
Let us note here for future use that, in view of Corollary \ref{corol}, we now have
\begin{align}
\label{p1}
\PP\left(\cN_o^+=1\right)&=\sum_{n\ge 1}n\pi_nz_\star^{n-1}(1-z_\star)=(1-z_\star){\varphi}'(z_\star)=(1-z_\star)^2\varphi'(1)\\
\label{wp1}
\widehat{\PP}\left(\cN_o^+=1\right)&=\sum_{n\ge 1}n\widehat{\pi}_nz_\star^{n-1}(1-z_\star)=(1-z_\star)\widehat{\varphi}'(z_\star)\\
\label{kappa}
\widehat{\EE}\left[\left.\cN_o^-\right|\cN_o^+=\cN_o^\star=0\right]&=\frac{1}{1-z_\star}\sum_{n=0}^\infty n\widehat{\pi}_nz_\star^n=\frac{z_\star\widehat{\varphi}'(z_\star)}{1-z_\star}.
\end{align}
Our first task consists in reducing the finiteness of $\EE\left[\beta^\star_o\right]$ to that of $\widehat{\EE}\left[\left.\frac{1}{\alpha_o}\right|\alpha_o>0\right]$. 
\begin{lem}[Reduction]Under (i), we have $\EE\left[\beta^\star_o\right]<\infty$ if and only if $\widehat{\EE}\left[\left.\frac{1}{\alpha_o}\right|\alpha_o>0\right]<\infty$.
\end{lem}
\begin{proof}[Proof of the ``only if'' part]
Since $\{\cN_o^+=1\}\subseteq\{\alpha_o=0\}\subseteq\{\beta_o^\star=\beta_o\}$, we have
\begin{align*}
\EE\left[\beta_o^\star\right]\ge\EE\left[\beta_o{\bf 1}_{(\cN_o^+=1)}\right]&=\EE\left[\frac{1}{\sum_{x\in\partial o}\alpha_{x\to o}}{\bf 1}_{(\cN_o^+=1)}\right]\\
&=\EE\left[\frac{1}{\sum_{x\in\partial o}\alpha_{x\to o}{\bf 1}_{(\alpha_{x\to o}>0)}}{\bf 1}_{(\cN_o^+=1)}\right]\\
&=\PP\left(\cN_o^+=1\right)\widehat{\EE}\left[\left.\frac{1}{\alpha_o}\right|\alpha_o>0\right].
\end{align*}
This is enough to conclude, since $\PP\left(\cN_o^+=1\right)>0$, by (\ref{p1}). 
\end{proof}
\begin{proof}[Proof of the ``if'' part]
Let us now assume that $\widehat{\EE}\left[\left.\frac{1}{\alpha_o}\right|\alpha_o>0\right]<\infty$, and verify that each term on the right-hand side of Formula (\ref{for:beta*}) is finite. For the first term, we write
\begin{align*}
\mathbf{E}\left[\mathbf{1}_{(\alpha_o=0)}\beta_o \right]
&=\mathbf{E}\left[\mathbf{1}_{(\beta_o<\infty)}\beta_o \right]\\
&=\mathbf{E}\left[\mathbf{1}_{(\cN_o^+\ge 1)}\frac{1}{\sum_{x\in\partial o}\alpha_{x\to o}} \right]
\end{align*}
where the first line follows from Lemma \ref{lm:i}, and the second from (\ref{decide:beta}) and (\ref{char:minus}).
Now, conditionally on $\cN_o^+$, the random variable $\sum_{x\in\partial o}\alpha_{x\to o}=\sum_{x\in\partial o}\alpha_{x\to o}{\bf 1}_{(\alpha_{x\to o}>0)}$ is distributed as the sum of $\cN_o^+$ i.i.d. random variables with law $\widehat{\PP}\left(\alpha_o\in\cdot|\alpha_o>0\right)$. Keeping only one of them yields
\begin{align}
\label{repeat}
\mathbf{E}\left[\mathbf{1}_{(\alpha_o=0)}\beta_o \right]&\le \PP\left(\cN_o^+\ge 1\right)\widehat{\EE}\left[\left.\frac{1}{\alpha_o}\right|\alpha_o>0\right],
\end{align}
which is finite. The second term is less than the first because $\{\cN_o^+\ge 2\}\subseteq\{\alpha_o=0\}$. For the third one, we observe the following: conditionally on $\cN_o^+$ and $\cN_o^-+\cN_o^\star$, the two random variables
\begin{align*}
\sum_{x \in \partial o} \alpha_{x \to o}\mathbf{1}_{(\alpha_{x \to o}>0)}\qquad\textrm{and}\qquad
\sum_{x \in \partial o} \beta_{x \to o}\mathbf{1}_{(\alpha_{x \to o}=0)},
\end{align*}
are independent, the first being distributed as a sum of $\cN_o^+$ i.i.d. random variables with law $\widehat{\PP}\left(\left.\alpha_o\in \cdot\right|\alpha_o>0\right)$, and the second as a sum of $\cN_o^\star+\cN_o^-$ i.i.d.  random variables with law $\widehat{\PP}\left(\left.\beta_o\in \cdot\right|\alpha_o=0\right)$. Keeping only one of the $\cN_o^+$ i.i.d. random variables in the first sum, we obtain
\begin{align*}
\EE\left[{\bf 1}_{(\cN_o^+\ge 2)}\frac{\sum_{x \in \partial o} \beta_{x \to o}\mathbf{1}_{(\alpha_{x \to o}=0)}}{\sum_{x \in \partial o} \alpha_{x \to o}\mathbf{1}_{(\alpha_{x \to o}>0)}}\right] \le \widehat{\EE}\left[\left.\frac{1}{\alpha_o}\right|\alpha_o>0\right]\EE\left[{\bf 1}_{(\cN_o^+\ge 2)}(\cN_o^-+\cN_o^\star)\right]\widehat{\EE}\left[\left.\beta_o\right|\alpha_o=0\right].
\end{align*}
The product on the right-hand side consists of three terms. The first is finite by assumption. The second is less than the expected degree at the root of our unimodular Galton-Watson tree, which is also finite. Finally, the inequality (\ref{repeat}) with $\PP$ replaced by $\widehat{\PP}$ shows that the third term is finite. 
\end{proof}
Since our running assumption (i) forces $M''(z_\star)\le 0$, the condition (ii) becomes $M''(z_\star)<0$ or equivalently,  $\widehat{\varphi}'(z_\star)<1$. To complete the proof of Theorem \ref{th:noext}, it therefore only remains to show that $\widehat{\EE}\left[\left.\frac{1}{\alpha_o}\right|\alpha_o>0\right]<\infty$ if and only if $\widehat{\varphi}'(z_\star)<1$, which we now do. 
\begin{lem} Under assumption (i), $\widehat{\EE}\left[\left.\frac{1}{\alpha_o}\right|\alpha_o>0\right]<\infty$ if and only if  $\widehat{\varphi}'(z_\star)<1$. 
\end{lem}
\begin{proof}[Proof of the ``only if'' part]On the one hand, using (\ref{decide:alpha}) and (\ref{char:plus}), we have
\begin{align*}
\widehat{\EE}\left[\left.\frac{1}{\alpha_o}\right|\alpha_o>0\right]&=1+\widehat{\EE}\left[\left.\sum_{x\in\partial o}\beta_{x\to o}\right|\alpha_o>0\right]\\
&=1+\widehat{\EE}\left[\left.\sum_{x\in\partial o}\beta_{x\to o}{\bf 1}_{(\beta_{x\to o}<\infty)}\right|\cN_o^+=\cN_o^\star=0\right]\\
&=1+\widehat{\EE}\left[\left.\cN_o^-\right|\cN_o^+=\cN_o^\star=0\right]\widehat{\EE}\left[\left.{\beta_o}\right|\beta_o<\infty\right]\\
&=1+\frac{\widehat{\varphi}'(z_\star)}{1-z_\star}\widehat{\EE}\left[{\beta_o}{\bf 1}_{(\beta_o<\infty)}\right],
\end{align*}
where the last line uses (\ref{kappa}) and (\ref{now}). On the other hand, using (\ref{decide:beta}), we have
\begin{align*}
\widehat{\EE}\left[\beta_o{\bf 1}_{(\cN_o^+=1)}\right]&=\widehat{\EE}\left[\frac{1}{\sum_{x\in\partial o}\alpha_{x\to o}}{\bf 1}_{(\cN_o^+=1)}\right]\\
&=\widehat{\EE}\left[\frac{1}{\sum_{x\in\partial o}\alpha_{x\to o}{\bf 1}_{(\alpha_{x\to o}>0)}}{\bf 1}_{(\cN_o^+=1)}\right]\\
&= \widehat{\PP}\left(\cN_o^+= 1\right)\widehat{\EE}\left[\left.\frac{1}{\alpha_o}\right|\alpha_o>0\right]\\
&= (1-z_\star)\widehat{\varphi}'(z_\star)\widehat{\EE}\left[\left.\frac{1}{\alpha_o}\right|\alpha_o>0\right],
\end{align*}
where the third line uses Corollary \ref{corol} and the last line uses (\ref{wp1}). 
Since $\{\cN_o^+=1\}\subseteq\{\beta_o<\infty\}$, we deduce from these two facts that
\begin{align*}
\widehat{\EE}\left[\left.\frac{1}{\alpha_o}\right|\alpha_o>0\right]\ge 1+\left(\widehat{\varphi}'(z_\star)\right)^2\widehat{\EE}\left[\left.\frac{1}{\alpha_o}\right|\alpha_o>0\right].
\end{align*}
The desired conclusion now clearly follows. 
\end{proof}
\begin{proof}[Proof of the ``if'' part]Fix $t>0$, and observe that by (\ref{decide}) and (\ref{char:plus}),
\begin{align}
\label{back0}
\widehat{\EE}\left[\left.\frac{1}{t\s_o(t)}\right|\alpha_o>0\right]&=\widehat{\EE}\left[\left.1+\sum_{x\in\partial o}\frac{\s_{x\to o}(t)}t\right|\alpha_o>0\right]\\
&=1+\widehat{\EE}\left[\left.\sum_{x\in\partial o}\frac{\s_{x\to o}(t)}t{\bf 1}_{(\beta_{x\to o}<\infty)}\right|\cN_o^+=\cN_o^\star=0\right]\\
&=1+\widehat{\EE}\left[\left.\cN_o^-\right|\cN_o^+=\cN_o^\star=0\right]\widehat{\EE}\left[\left.\frac{\s_o(t)}{t}\right|\beta_o<\infty\right]\\
&= 1+\frac{z_\star\widehat{\varphi}'(z_\star)}{1-z_\star}\widehat{\EE}\left[\left.\frac{\s_o(t)}{t}\right|\beta_o<\infty\right]
\end{align}
where the last line uses (\ref{kappa}). On the other hand, using (\ref{decide}) and (\ref{char:minus}), we have
\begin{align}
\nonumber
\widehat{\EE}\left[\left.\frac{\s_o(t)}{t}\right|\beta_o<\infty\right]&=\widehat{\EE}\left[\left.\left(\sum_{x\in\partial o}t\s_{x\to o}(t)\right)^{-1}\,\right|\cN_o^+\ge 1\right]
\\ 
\label{back1}
&\le \widehat{\EE}\left[\left.\left(\sum_{x\in\partial o}t\s_{x\to o}(t){\bf 1}_{(\alpha_{x\to o}>0)}\right)^{-1}\,\right|\cN_o^+\ge 1\right].
% = \int_0^1\frac{1}{z}\phi\left(\widehat{\EE}\left[\left.z^{t\s_o(t)}\right|\alpha_o>0\right]\right)dz,
\end{align}
Conditionally on $\cN_o^+$, the integrand on the right-hand side is distributed as the reciprocal of the sum of $\cN_o^+$ i.i.d. random variables with law $\widehat{\PP}\left(\left.t\s_o(t)\in\cdot \right|\alpha_o>0\right)$. To exploit this i.i.d. structure, we transform the reciprocal $(\cdot)^{-1}$ into a power via the trivial identity
\begin{align*}
r^{-1} = \int_0^1 z^{r-1}\,{\rm d} z,
\end{align*}
valid for any $r>0$. With $r=\sum_{x\in\partial o}t\s_{x\to o}(t){\bf 1}_{(\alpha_{x\to o}>0)}$, we obtain
\begin{align*}
\widehat{\EE}\left[\left.\left(\sum_{x\in\partial o}t\s_{x\to o}(t){\bf 1}_{(\alpha_{x\to o}>0)}\right)^{-1}\,\right|\cN_o^+\ge 1\right]&=\int_0^1\widehat{\EE}\left[\left.z^{\sum_{x\in\partial o}t\s_{x\to o}(t){\bf 1}_{(\alpha_{x\to o}>0)}}\right|\cN_o^+\ge 1\right]\,\frac{{\rm d} z}{z}\\
&= \sum_{n=1}^\infty\widehat{\PP}\left(\cN^+_o=n|\cN_o^+\ge 1\right)\int_0^1\left(\widehat{\EE}\left[\left.z^{t\s_{o}(t)}\right|\alpha_{o}>0\right]\right)^n\frac{{\rm d} z}{z}.
\end{align*}
We now fix some $\varepsilon\in(0,1)$ and $n\ge 1$, and estimate the integral on the right-hand side by splitting it into two parts: for $z\in (\varepsilon,1)$, we use the crude bound $\widehat{\EE}\left[\left.z^{t\s_{o}(t)}\right|\alpha_{o}>0\right]\le 1$ to obtain
\begin{align*}
\int_\varepsilon^1\left(\widehat{\EE}\left[\left.z^{t\s_{o}(t)}\right|\alpha_{o}>0\right]\right)^n\frac{{\rm d} z}{z}\le \ln\left(\frac{1}{\varepsilon}\right).
%+\widehat{\EE}\left[\left.\varepsilon^{\alpha_{o}}\right|\alpha_{o}>0\right]^{n-1}\int_0^1 \widehat{\EE}\left[\left.z^{t\s_{o}(t)}\right|\alpha_{o}>0\right]\frac{{\rm d} z}{z},
\end{align*}
For $z\in(0,\varepsilon)$, we use the observation that  $z^{t\s_o(t)}\le \varepsilon^{\alpha_o}$ to write
\begin{align*}
\int_0^\varepsilon\left(\widehat{\EE}\left[\left.z^{t\s_{o}(t)}\right|\alpha_{o}>0\right]\right)^n\frac{{\rm d} z}{z}
&\le \widehat{\EE}\left[\left.\varepsilon^{\alpha_{o}}\right|\alpha_{o}>0\right]^{n-1}\int_0^\varepsilon \widehat{\EE}\left[\left.z^{t\s_{o}(t)}\right|\alpha_{o}>0\right]\frac{{\rm d} z}{z}\\
&\le \widehat{\EE}\left[\left.\varepsilon^{\alpha_{o}}\right|\alpha_{o}>0\right]^{n-1}\widehat{\EE}\left[\left.\frac{1}{t\s_{o}(t)}\right|\alpha_{o}>0\right].
\end{align*}
Inserting these estimates into the above series and recalling (\ref{back0})-(\ref{back1}), we arrive at
\begin{align}
\label{conclude}
\widehat{\EE}\left[\left.\frac{1}{t\s_o(t)}\right|\alpha_o>0\right]\le 1+\frac{z_\star\widehat{\varphi}'(z_\star)}{1-z_\star}\left(\log\left(\frac{1}{\varepsilon}\right)+ \Phi\left(\widehat{\EE}\left[\left.\varepsilon^{\alpha_o}\right|\alpha_o>0\right]\right) \widehat{\EE}\left[\left.\frac{1}{t\s_o(t)}\right|\alpha_o>0\right]\right),
\end{align}
where we have introduced the short-hand
\begin{align*}
\Phi(u)=\widehat{\EE}\left[\left.u^{\cN^+_o-1}\right|\cN^+_o \ge 1\right]=\sum_{n=0}^\infty\widehat{\PP}\left(\cN^+_o=n+1|\cN_o^+\ge 1\right)u^n.
\end{align*}
Now, observe that
\begin{align*}
\Phi(0)=\frac{\widehat{\PP}\left(\cN^+_o=1\right)}{\widehat{\PP}\left(\cN_o^+\ge 1\right)}=\frac{1-z_\star}{z_\star}\widehat{\varphi}'(z_\star),
\end{align*}
where we have used (\ref{char:minus}), (\ref{now}) and (\ref{wp1}).
By continuity of $\Phi$, we deduce that
\begin{align*}
\frac{z_\star\widehat{\varphi}'(z_\star)}{1-z_\star}\Phi\left(\widehat{\EE}\left[\left.\varepsilon^{\alpha_o}\right|\alpha_o>0\right]\right)\xrightarrow[\varepsilon\to 0]{}\left(\widehat{\varphi}'(z_\star)\right)^2.
\end{align*}
If $\widehat{\varphi}'(z_\star)<1$, we can  choose $\varepsilon>0$ so that $\frac{z_\star\widehat{\varphi}'(z_\star)}{1-z_\star}\Phi\left(\widehat{\EE}\left[\left.\varepsilon^{\alpha_o}\right|\alpha_o>0\right]\right)<1$ and rewrite (\ref{conclude}) as
\begin{align*}
\widehat{\EE}\left[\left.\frac{1}{t\s_o(t)}\right|\alpha_o>0\right]\le \frac{1+\frac{z_\star\widehat{\varphi}'(z_\star)}{1-z_\star} \log\left(\frac{1}{\varepsilon}\right)}{1-\frac{z_\star\widehat{\varphi}'(z_\star)}{1-z_\star}\Phi\left(\widehat{\EE}\left[\left.\varepsilon^{\alpha_o}\right|\alpha_o>0\right]\right)}.
\end{align*}
The right-hand side is finite and independent of $t$, so letting $t\to 0$ concludes the proof.
\end{proof}

\subsection{Proof of Theorem \ref{th:ext}}
We use the following Lemma, whose proof is trivial once we observe that a sequence of non-negative random variables $(X_n)_{n\ge 1}$ tends to $\infty$ in probability if and only if $\EE[\exp(-X_n)]\xrightarrow[n\to\infty]{} 0$.
\begin{lem}\label{lm:proba}
Let $k\ge 1$ be a fixed integer,  and let $(X^{(1)}_n)_{n \geqslant 1},\ldots,(X^{(k)}_n)_{n \geqslant 1}$ be $k$ i.i.d. copies of an arbitrary sequence  $(X_n)_{n \geqslant 1}$  of non-negative random variables. Then the sequence $(Y_n)_{n \geqslant 1}$ defined by $Y_n:=X^{(1)}_n + \dotsb + X^{(k)}_n$ tends to $\infty$ in probability if and only if $(X_n)_{n \geqslant 1}$ does.
\end{lem}

\begin{proof}[Proof of Theorem \ref{th:ext}]
Assume that condition (i) fails. By Lemma \ref{lm:i}, this means that the event
\begin{align*}
\cE:=\left\{\alpha_o=0,\beta_o=\infty\right\}=\{\cN_o^+=0,\cN_o^\star\ge 1\}
\end{align*}
has positive probability under $\PP$ and $\widehat{\PP}$.  On this event, the recursion (\ref{decide}) can be rewritten as
 \begin{align*}
\s^\star_o(t)=
\left(t+\sum_{x\in\partial o}\s_{x\to o}(t){\bf 1}_{(\beta_{x\to o}<\infty)}+\sum_{x\in\partial o}\s_{x\to o}^\star(t){\bf 1}_{(\alpha_{x\to o}=0,\beta_{x\to o}=\infty)}\right)^{-1}.
\end{align*}
The first sum on the right-hand side tends to $0$ as $t\to 0$ by definition of $\beta_{x\to o}$. On the other hand, by Corollary \ref{corol}, conditionally on $\cN_o^\star$, the second sum is distributed as the sum of $\cN_o^\star$ i.i.d. variables with law $\widehat{\PP}\left(\s_o^\star\in\cdot|\cE\right)$. We emphasize that this statement is valid under both $\PP$ and $\widehat{\PP}$ (only the distribution of $\cN_o^\star$ differ). Applying Lemma \ref{lm:proba} to both situations, we deduce that along any deterministic sequence $(t_n)_{n\ge 1}$ of positive numbers with $t_n\to 0$ as $n\to\infty$, the following conditions are equivalent:
\begin{enumerate}
\item[(a)] $\s_o^\star(t_n)\xrightarrow[n\to\infty]{}0$ in probability under $\PP(\cdot|\cE)$;
\item[(b)] $\s_o^\star(t_n)\xrightarrow[n\to\infty]{}\infty$
in probability under $\widehat{\PP}\left(\cdot|\cE\right)$;
\item[(c)] $\s_o^\star(t_n)\xrightarrow[n\to\infty]{}0$ in probability under $\widehat{\PP}\left(\cdot|\cE\right)$. 
\end{enumerate}
Of course, (b) and (c) are incompatible, and so (a) can never hold. In particular, this rules out the possibility that $\EE[\s^\star_o(t_n)]\to 0$ as $n\to\infty$, and since $(t_n)_{n\ge 1}$ is arbitrary, 
we conclude that
\begin{align}
\label{lastgoal}
\liminf_{t\to 0^+}\EE\left[\s^\star_o(t)\right]>0.
\end{align}
By Lemma \ref{lm:ext}, this is more than enough to ensure that $\mu_{\UGWT(\pi)}$ has extended states at zero. 
\end{proof}
\bibliography{bibli}

\end{document}